\documentclass[11pt]{article}
\usepackage{graphicx}
\usepackage{amsmath,amsfonts,amssymb}
\usepackage{citesort}
\usepackage{url}
\usepackage{amsmath}
\usepackage{graphicx}
\usepackage{epsfig}
\usepackage{color}
\def\tto{\;{\lower 1pt \hbox{$\rightarrow$}}\kern -10pt
\hbox{\raise 2pt \hbox{$\rightarrow$}}\;}
\def\Hat{\widehat}

\def\ra{\rangle}
\def\la{\langle}
\def\ve{\varepsilon}
\def\B{I\!\!B}
\def\h{\hfill\Box}
\def\R{\Bbb R}
\def\N{I\!\!N}
\def\ox{\bar{x}}

\def\oy{\bar{y}}
\def\oz{\bar{z}}

\def\gph{\mbox{\rm gph}\,}
\def\epi{\mbox{\rm epi}\,}

\def\h{\hfill\square}
\def\dn{\downarrow}
\def\O{\Omega}
\def\ph{\varphi}

\def \N{I\!\!N}

\newcounter{lk}

\begin{document}
\begin{center}
{\bf LIPSCHITZ PROPERTIES OF NONSMOOTH FUNCTIONS AND SET-VALUED MAPPINGS VIA GENERALIZED DIFFERENTIATION AND APPLICATIONS}\\[2ex]
Nguyen Mau Nam\footnote{Fariborz Maseeh Department of
Mathematics and Statistics, Portland State University, Portland, OR 97202, United States (mau.nam.nguyen@pdx.edu). The research of Nguyen Mau Nam was partially supported by
the Simons Foundation under grant \#208785.} and Gerardo Lafferriere\footnote{Fariborz Maseeh Department of
Mathematics and Statistics, Portland State University, Portland, OR 97202, United States (gerardoL@pdx.edu).}
\end{center}
\small{\bf Abstract.} In this paper, we revisit the Mordukhovich's subdifferential criterion for Lipschitz continuity of nonsmooth functions and coderivative criterion for the Aubin/Lipschitz-like property of set-valued mappings in finite dimensions. The criteria are useful and beautiful results in modern variational analysis showing the state of the art of the field. As an application, we establish necessary and sufficient conditions for Lipschitz continuity of the minimal time function and the scalarization function, that play an important role in many aspects of nonsmooth analysis and optimization. \\[1ex]
{\bf Key words.} Lipschitz property, generalized differentiation, subdifferential, coderivative, minimal time function, scalarization function\\[1ex]
{\bf AMS subject classifications.} 49J53, 49J52, 90C31

\newtheorem{Theorem}{Theorem}[section]
\newtheorem{Proposition}[Theorem]{Proposition}
\newtheorem{Remark}[Theorem]{Remark}
\newtheorem{Lemma}[Theorem]{Lemma}
\newtheorem{Corollary}[Theorem]{Corollary}
\newtheorem{Definition}[Theorem]{Definition}
\newtheorem{Example}[Theorem]{Example}
\renewcommand{\theequation}{\thesection.\arabic{equation}}
\normalsize

\section{Introduction and Preliminaries}

Lipschitz continuity is an important concept in mathematical analysis. In modern variational analysis, it has been generalized for set-valued mappings. Among many extensions, the \emph{pseudo-Lipschitzian property} introduced by Aubin \cite{af} has been well recognized as a natural and useful one. It is now called by different names such as the \emph{Aubin property} or the \emph{Lipschitz-like property}. The concept has been used extensively in the study of sensitivity analysis of optimization problems and variational inequalities. It also plays an important role on developing generalized differentiation calculi for nonsmooth functions and set-valued mappings; see \cite{mor,rw} and the references therein for more discussions on the history of the concept, as well as many important applications to variational analysis, optimization, and optimal control. \vspace*{0.05in}

Recall that a set-valued mapping $\mathcal{F}: {\Bbb R}^m\tto {\Bbb R}^n$ has the \emph{Aubin property} around  $(\ox,\oy)\in \gph \mathcal{F}:=\{(x,y)\in \R^m\times\R^n\; |\; y\in \mathcal{F}(x)\}$ if there exist neighborhoods $V$ of $\ox$, $W$ of $\oy$, and a constant $\ell\geq 0$ such that
\begin{equation*}
\mathcal{F}(x)\cap W\subseteq \mathcal{F}(u)+\ell \|x-u\|\B \mbox{ for all }x, u\in V.
\end{equation*}
The first effort to characterize the Aubin property using generalized differentiation was made by Rockafellar \cite{r85}. A sufficient condition for $\mathcal{F}$ to have the Aubin property around $(\ox,\oy)\in \gph \mathcal{F}$ was given as follows:
\begin{equation*}
[(u,0)\in N_C((\ox,\oy);\gph \mathcal{F})]\Rightarrow u=0,
\end{equation*}
where $N_C((\ox,\oy);\gph \mathcal{F})$ is the \emph{Clarke normal cone} to $\gph \mathcal{F}$ at $(\ox,\oy)$; see \cite{c-1983}.

However, the Clarke normal cone is \emph{too large} to be able to recognize the Aubin property of set-valued mappings in many different settings; see more discussions in \cite{m92,m93}. Therefore, it was a need to find a necessary and sufficient condition for this property using \emph{smaller normal cone structures}, and the \emph{Mordukhovich/limiting normal cone} \cite{m76} gives an answer to this question. The implication
\begin{equation}\label{c1}
[(u,0)\in N((\ox,\oy);\gph \mathcal{F})]\Rightarrow u=0
\end{equation}
in terms of the Mordukhovich normal cone $N((\ox,\oy);\gph \mathcal{F})$ to $\gph \mathcal{F}$ at $(\ox,\oy)$ is indeed a necessary and sufficient condition for $\mathcal{F}$ to have the Aubin property around $(\ox,\oy)$. This striking result was first proved by Mordukhovich in \cite[Theorem 5.4]{m92}. It is now called the\emph{ Mordukhovich's coderivative criterion} for the Aubin property. We will get back to the idea behind Mordukhobvich's proof after presenting some important concepts of variational analysis. The readers are referred to \cite{mor,rw} for more detail.\vspace*{0.05in}

Let $\mathcal{F}: {\Bbb R}^m\tto {\Bbb R}^n$ be a set-valued mapping. The \emph{domain} of the mapping is 
\begin{equation*}
\mbox{dom }\mathcal{F}:=\{x\in {\Bbb R}^m\; |\; \mathcal{F}(x)\neq \emptyset\}.
\end{equation*}
Given $x\in \R^n$ and a subset $\O\subseteq {\Bbb R}^n$, the \emph{distance function} from $x$ to $\O$ is given by
\begin{equation*}
d(x;\O):=\inf\{\|x-\omega\|\;|\; \omega\in \O\}.
\end{equation*}
The set
\begin{equation*}
\Pi(x;\O):=\{\omega\in\O\; |\; d(x;\O)=\|x-\omega\|\}
\end{equation*}
is called the \emph{metric projection} from $x$ to $\O$.

The following function defined on ${\Bbb R}^m\times {\Bbb R}^n$ will play an important role throughout the paper:
\begin{equation}\label{scalar}
D(x, y):=d(y; \mathcal{F}(x)).
\end{equation}
Let $\O\subseteq {\Bbb R}^n$ and let $\ox\in \O$. A vector $v\in {\Bbb R}^n$ is called a \emph{Fr\'echet normal} to $\O$ at $\ox$ if
\begin{equation*}
\la v, x-\ox\ra \leq \mbox{\rm o}(\|x-\ox\|) \mbox{ for }x\in \O.
\end{equation*}
The set of all Fr\'echet normals to $\O$ at $\ox$ is called the \emph{Fr\'echet normal cone} to $\O$ at $\ox$, denoted by $\Hat N(\ox;\O)$. \vspace*{0.05in}

A vector $v\in {\Bbb R}^n$ is called a \emph{limiting normal} to $\O$ at $\ox$ if there are sequences $x_k\xrightarrow{\O}\ox$ and $v_k\rightarrow v$ with $v_k\in \Hat N(x_k;\O)$. In this definition, the notation $x_k\xrightarrow{\O}\ox$ means that $x_k\to\ox$ and $x_k\in \O$ for every $k$. The set of all limiting normals to $\O$ at $\ox$ is called the \emph{Mordukhovich/limiting normal cone} to the set at $\ox$ and is denoted by $N(\ox; \O)$.  \vspace*{0.05in}

Let $\psi: {\Bbb R}^n\to (-\infty, \infty]$ be an extended real-valued function and let $\ox$ be an element of the \emph{domain} of the function $\mbox{\rm dom }\psi:=\{x\in {\Bbb R}^n\; |\; \psi(x)<\infty\}$. The \emph{Fr\'echet subdifferential} of $\psi$ at $\ox$ is defined by
\begin{equation*}
\Hat\partial\psi(\ox):=\{v\in {\Bbb R}^n\; |\; \la v, x-\ox\ra \leq  \psi(x)-\psi(\ox)+\mbox{\rm o}(\|x-\ox\|)\}.
\end{equation*}
The \emph{limiting/Mordukhovich subdifferential} of $\psi$ at $\ox$, denoted by $\partial \psi(\ox)$, is the set of vectors $v\in {\Bbb R}^n$ such that there exist sequences $x_k\xrightarrow{\psi}\ox$, and $v_k\in \Hat\partial\psi(x_k)$ with $v_k\rightarrow v$. The \emph{singular subdifferential} of $\psi$ at $\ox$, denoted by $\partial^\infty\psi(\ox)$, is the set of all vectors $v\in {\Bbb R}^n$ such that there exist sequences $\lambda_k\dn 0$, $x_k\xrightarrow{\psi}\ox$, and $v_k\in \Hat\partial\psi(x_k)$ with $\lambda_kv_k\rightarrow v$. In these definitions, $x_k\xrightarrow{\psi}\ox$ means that $x_k\to\ox$ and $\psi(x_k)\to\psi(\ox)$, and $\lambda_k\dn 0$ means that $\lambda_k\to 0$ and $\lambda_k\geq 0$ for every $k$. Both Fr\'echet and limiting subdifferential constructions reduce to the subdifferential in the sense of convex analysis when the function involved is convex. Moreover, if $\psi$ is lower semicontinuous around $\ox$, one has the following representation:
\begin{equation*}
\partial^\infty\psi(\ox)=\{v\in \R^n\; |\; (v,0)\in N((\ox, \psi(\ox)); \epi \psi)\}.
\end{equation*}

An extended real-valued function $\psi: {\Bbb R}^n\to (-\infty, \infty]$ is called \emph{Lipschitz continuous} around $\ox\in \mbox{\rm dom }\psi$ if there exist a constant $\ell$ and a neighborhood $V$ of $\ox$ such that
\begin{equation*}
|\psi(x)-\psi(u)|\leq \ell \|x-u\| \mbox{ for all }x,u\in V.
\end{equation*}
If this equality is replaced by
\begin{equation*}
|\psi(x)-\psi(\ox)|\leq \ell \|x-\ox\| \mbox{ for all }x\in V,
\end{equation*}
we say that $\psi$ is \emph{calm} at $\ox$.

Let $\mathcal{F}: {\Bbb R^m}\tto {\Bbb R}^n$ be a set-valued mapping and let $(\ox,\oy)\in
\gph \mathcal{F}$. The \emph{Mordukhovich/limiting coderivative} of $\mathcal{F}$ at $(\ox,\oy)$ is a set-valued
mapping, denoted by $D^*\mathcal{F}(\ox,\oy): {\Bbb R}^n\tto {\Bbb R}^m$ defined by
\begin{equation*}
D^*\mathcal{F}(\ox,\oy)(v):=\{u\in {\Bbb R}^m \; |\; (u, -v)\in N((\ox,\oy);
\gph \mathcal{F})\}.
\end{equation*}
The necessary and sufficient condition (\ref{c1}) for $\mathcal{F}$ to have the Aubin property around $(\ox,\oy)\in \gph \mathcal{F}$ can be equivalently represented in terms of the mordukhovich coderivative in the theorem below.
\begin{Theorem}\label{mcriterion} {\rm (\cite[Theorem 5.4]{m92})} Let $\mathcal{F}: {\Bbb R}^m\tto {\Bbb R}^n$ be a set-valued mapping whose graph is locally closed around $(\ox,\oy)\in \gph \mathcal{F}$. Then $\mathcal{F}$ has the Aubin property around $(\ox,\oy)$ if and only if $$D^*\mathcal{F}(\ox,\oy)(0)=\{0\}.$$
\end{Theorem}

The necessary part of the proof of Theorem \ref{mcriterion} by Mordukhovich is based on a direct estimate of $\|u\|$ in terms of $\|v\|$ for $u\in D^*\mathcal{F}(\ox,\oy)(v)$ from \cite[Propostion 5.2]{m92}. The sufficient part of the proof of the theorem is based on an upper estimate for subgradients of a general marginal function specified to the case of function (\ref{scalar}).  Other important tools for the proof are the necessary and sufficient condition for the Lipschitz continuity of extended real-valued functions, and the equivalence between the Aubin property of $\mathcal{F}$ at $(\ox,\oy)$ and the Lipschitz continuity of function $D$ given by {\rm (\ref{scalar})} at the same point in the theorems below.

\begin{Theorem}\label{scarl}{\rm (\cite[Theorem 3.2]{r85})} Let $\mathcal{F}: {\Bbb R}^m\tto {\Bbb R}^n$ be a set-valued mapping whose graph is locally closed around $(\ox,\oy)\in \gph \mathcal{F}$. Then $\mathcal{F}$ has the Aubin property around $(\ox,\oy)$ if and only if the function $D$ given by {\rm (\ref{scalar})} is Lipschitz continuous around $(\ox,\oy)$.
\end{Theorem}

\begin{Theorem}\label{mcriterionf}{\rm \cite[Theorem 2.1]{m88})} Let $\psi: {\Bbb R}^n\to (-\infty, \infty]$ be an extended real-valued function that is lower semicontinuous around $\ox\in \mbox{\rm dom }\psi$. Then $\psi$ is Lipschitz continuous around $\ox$ if and only if $$\partial^\infty\psi(\ox)=\{0\}.$$
\end{Theorem}

A self-contained proof was given by Rockafellar and Wets \cite[Theorem 9.40]{rw}. However, their proof is not easy to understand, especially for those who have just started the study of variational analysis. Thus, our first goal in this paper is to provide a simpler self-contained proof of the Mordukhovich coderivative criterion. To achieve the goal, in Section 2 of the paper, we will show that the subdifferential criterion for Lipschitz continuity of nonsmooth functions, Theorem \ref{mcriterionf}, and the coderivative criterion for the Aubin property of set-valued mappings, Theorem \ref{mcriterion}, are in fact equivalent. Instead of using \cite[Theorem 4.1]{m92} as in \cite[Theorem 5.4]{m92} , we provide a simple direct way to obtain the upper estimate for singular subgradients of the distance function (\ref{scalar}) in terms of $D^*\mathcal{F}(\ox,\oy)(0)$. Then we will show that the estimate becomes an equality under the Lipschitz continuity of function (\ref{scalar}), and the equality will be used to obtain the necessary condition. As an application, in Section 3, we use these criteria to study Lipschitz continuity of the \emph{minimal time function} and the \emph{scalarization function}, that play a crucial role in many aspects of nonsmooth analysis and optimization; see, e.g., \cite{cowo,MN} and the references therein.\vspace*{0.05in}

Throughout the paper, $\la \cdot, \cdot, \ra$ denotes the dot product in $\R^n$; $\B$ stands for the closed unit ball of $\R^n$; $\B(\ox; r)$ denotes the closed ball with center at $\ox\in \R^n$ and radius $r$. We will use the ``sum'' norm in ${\Bbb R}^m\times{\Bbb R}^n$, and use the Euclidean norm in ${\Bbb R}^m$ and ${\Bbb R}^n$.

\section{Lipschitz Properties via Generalized Differentiation}

Let us start with a simple proof of a known result. For simplicity, we assume the closedness of the graph of the set-valued mapping $\mathcal{F}$ instead of the local closedness.

\begin{Proposition}\label{cg} Let $\mathcal{F}: {\Bbb R}^m\tto {\Bbb R}^n$ be a set-valued mapping with closed graph. Then the function $D$ given by {\rm (\ref{scalar})} is lower semicontinuous.
\end{Proposition}
{\bf Proof.} For any $\alpha\in \R$, one has the following representation of the \emph{$\alpha-$level set}:
\begin{align*}
\{(x,y)\in {\Bbb R}^m\times {\Bbb R}^n\; |\; D(x,y)\leq \alpha\}&=\{(x,y)\in {\Bbb R}^m\times {\Bbb R}^n\; |\; d(y; \mathcal{F}(x))\leq \alpha\}\\
&=\{(x,y)\in {\Bbb R}^m\times {\Bbb R}^n\; |\; y\in \mathcal{F}(x)+\alpha\B\}\\
&=\gph \mathcal{G},
\end{align*}
where $\mathcal{G}(x):=\mathcal{F}(x)+\alpha\B$. Since $\gph \mathcal{F}$ is closed, $\gph \mathcal{G}$ is closed. Indeed, let $(x_k, y_k)\xrightarrow{\gph G}(\ox,\oy)$. Then $y_k=v_k+\alpha e_k$, where $v_k\in \mathcal{F}(x_k)$ and $\|e_k\|\leq 1$. Then $(v_k)$ and $(e_k)$ are bounded, and we can assume, without loss of generality, that $v_k\to \bar v$ and $e_k\to\bar e\in \B$. Thus, $\oy=\bar v+\alpha\bar e\in F(\ox)+\alpha\B$ as $\gph \mathcal{F}$ is closed. Since any $\alpha-$level set is closed, $D$ is lower semicontinuous. $\h$

\begin{Proposition}\label{l1} Consider the function $D$ given by {\rm (\ref{scalar})}, where $\mathcal{F}: {\Bbb R}^m\tto {\Bbb R}^n$ is a set-valued mapping with closed graph. Then
\begin{equation*}
[(v, w)\in \Hat\partial D(\ox,\oy)]\Rightarrow [(v,w)\in \Hat N((\ox,\oz); \gph \mathcal{F}),\; \|w\|\leq 1]
\end{equation*}
for any $\oz\in \Pi(\oy; \mathcal{F}(\ox))$.
\end{Proposition}
{\bf Proof.} Fix $(v, w)\in \Hat\partial D(\ox,\oy)$ and $\oz\in \Pi(\oy, \mathcal{F}(\ox))$. Then $d(\oy; \mathcal{F}(\ox))=\|\oy-\oz\|$, and for any $\ve>0$, there exists $\delta>0$ such that
\begin{equation}\label{1}
\la v, x-\ox\ra +\la w, y-\oy\ra \leq d(y; \mathcal{F}(x))-d(\oy; \mathcal{F}(\ox))+\ve (\|x-\ox\|+\|y-\oy\|)
\end{equation}
whenever $\|x-\ox\|<\delta$ and $\|y-\oy\|<\delta$.

Fix any $(x,y)\in \gph \mathcal{F}$ with $\|x-\ox\|<\delta$ and $\|y-\oz\|<\delta$. Then $\|(y-\oz+\oy)-\oy\|<\delta$. Thus, by (\ref{1}) with $y$ being replaced by $y-\oz+\oy$,
\begin{align*}
\la v, x-\ox\ra +\la w, y-\oz\ra &\leq d(y-\oz+\oy; \mathcal{F}(x))-d(\oy; \mathcal{F}(\ox))+\ve (\|x-\ox\|+\|y-\oz\|)\\
&\leq d(y; \mathcal{F}(x))+\|\oz-\oy\|-\|\oy-\oz\|+\ve (\|x-\ox\|+\|y-\oz\|)\\
&=\ve (\|x-\ox\|+\|y-\oz\|),
\end{align*}
since $d(y; \mathcal{F}(x))=0$ and $d(\oy; \mathcal{F}(\ox))=\|\oz-\oy\|$. It follows that $(v, w)\in \Hat N((\ox,\oz); \gph \mathcal{F})$. Using $x\equiv\ox$ in (\ref{1}), one one has $w\in\Hat\partial d(\oy;\O)$, where $\O:=\mathcal{F}(\ox)$. Since $d(\cdot;\O)$ is Lipschitz continuous with constant $\ell=1$, one sees that $\|w\|\leq 1$. The proof is now complete. $\h$ \vspace*{0.05in}

In Proposition \ref{l1}, we can replace the condition $\|w\|\leq 1$ by $\|w\|=1$ if $\oy\notin \mathcal{F}(\ox)$, but we do not need this in our subsequent analysis. Let us continue with another useful result for proving the necessary condition of the Mordukhovich's coderivative criterion. In the proposition below, we will assume the Lipschitz continuity of function (\ref{scalar}) although it is possible to make this assumption weaker.

\begin{Proposition}\label{l3} Let $\mathcal{F}: \Bbb R^m\tto \Bbb R^n$ be a set-valued mapping with closed graph. Suppose that the function $D$ given by {\rm (\ref{scalar})} is Lipschitz continuous around $(\ox,\oy)\in \gph \mathcal{F}$. Then
\begin{equation*}
[(v, w)\in \Hat N((\ox,\oy); \gph \mathcal{F})]\Rightarrow [(v,w)\in \lambda \Hat\partial D(\ox,\oy), \mbox{\rm where } \lambda:=\|w\|].
\end{equation*}
\end{Proposition}
{\bf Proof.} Fix any $(v, w)\in \Hat N((\ox,\oy); \gph \mathcal{F})$. Then for any $\ve>0$, there exists $\delta>0$ such that
\begin{equation*}
\la v, x-\ox\ra +\la w, z-\oy\ra \leq \ve (\|x-\ox\|+\|z-\oy\|)
\end{equation*}
whenever $(x,z)\in \gph \mathcal{F}$, $\|x-\ox\|<\delta$, $\|z-\oy\|<\delta$.

Fix any sequence $(x_k, y_k)$ that converges to $(\ox, \oy)$, $(x_k, y_k)\neq (\ox,\oy)$. Since $D$ is Lipschitz continuous around $(\ox, \oy)$, it is finite around this point, and hence $\mathcal{F}(x_k)$ is nonempty (and closed) around $\ox$. Pick $z_k\in \Pi(y_k; \mathcal{F}(x_k))$. Then $(x_k,z_k)\in \gph \mathcal{F}$, and
\begin{equation*}
D(x_k,y_k)=d(y_k; \mathcal{F}(x_k))=\|y_k-z_k\|.
\end{equation*}
Since $D$ is Lipschitz continuous around $(\ox,\oy)$ and $D(\ox,\oy)=d(\oy; \mathcal{F}(\ox))=0$, one has $z_k\to\oy$.  Let $\ell$ be a Lipschitz constant of $D$ around $(\ox,\oy)$. For sufficiently large $k$, the following holds
$$d(y_k; \mathcal{F}(x_k))=d(y_k; \mathcal{F}(x_k))-d(\oy; \mathcal{F}(\ox))\leq \ell (\|x_k-\ox\|+\|y_k-\oy\|).$$
It follows that
\begin{align*}
\la v, x_k-\ox\ra +\la w, y_k-\oy\ra & = \la v, x_k-\ox\ra +\la w, z_k-\oy\ra+ \la w, y_k-z_k\ra\\
&\leq \ve (\|x_k-\ox\|+\|z_k-\oy\|)+\la w, y_k-z_k\ra\\
&\leq \ve(\|x_k-\ox\|+\|y_k-\oy\|+\|z_k-y_k\|)+\|w\|\; \|y_k-z_k\|\\
&\leq \ve (\|x_k-\ox\|+\|y_k-\oy\|)+(\ve +\|w\|)\|y_k-z_k\|\\
&= \ve (\|x_k-\ox\|+\|y_k-\oy\|)+(\ve +\|w\|)d(y_k; \mathcal{F}(x_k))\\
&= \ve (\|x_k-\ox\|+\|y_k-\oy\|) +\ve d(\mathcal{F}(x_k), y_k) + \|w\|d(y_k; \mathcal{F}(x_k))\\
&\leq \ve(\ell+1) (\|x_k-\ox\|+\|y_k-\oy\|) +\lambda d(y_k; \mathcal{F}(x_k))\\
&= \ve(\ell+1) (\|x_k-\ox\|+\|y_k-\oy\|) +\lambda d(y_k; \mathcal{F}(x_k))-\lambda d(\oy, \mathcal{F}(\ox))\\
&= \ve(\ell+1)( \|x_k-\ox\|+\|y_k-\oy\|) +\lambda D(x_k,y_k)-\lambda D(\ox,\oy).
\end{align*}
Since $\ve>0$ is arbitrary, one has
\begin{equation*}
\liminf_{k\to\infty} \dfrac{\lambda D(x_k,y_k)-\lambda D(\ox,\oy)-\la v, x_k-\ox\ra -\la w, y_k-\oy\ra}{\|x_k-\ox\|+\|y_k-\oy\|}\geq 0.
\end{equation*}
Thus,
\begin{equation*}
\liminf_{(x,y)\to (\ox,\oy)} \dfrac{\lambda D(x,y)-\lambda D(\ox,\oy)-\la v, x-\ox\ra -\la w, y-\oy\ra}{\|x-\ox\|+\|y-\oy\|}\geq 0.
\end{equation*}
Therefore, $(v,w)\in \lambda \Hat\partial D(\ox,\oy), \mbox{\rm where } \lambda=\|w\|$. $\h$

\begin{Proposition}\label{l2} Let $\mathcal{F}$ be a set-valued mapping with closed graph. For any $(\ox,\oy)\in \gph \mathcal{F}$, one has
\begin{equation}\label{upper estimate}
\partial^\infty D(\ox,\oy)\subseteq \{(x^*,0)\; |\; x^*\in D^*\mathcal{F}(\ox,\oy)(0)\}.
\end{equation}
The equality holds if $D$ is Lipschitz continuous around $(\ox,\oy)$.
\end{Proposition}
{\bf Proof.} Fix any $(v, w)\in \partial^\infty D(\ox,\oy)$. Then there exist sequences $\lambda_k\dn 0$, $(x_k, y_k)\xrightarrow{D}(\ox,\oy)$, $(v_k, w_k)\in \Hat\partial D(x_k,y_k)$ with
\begin{equation*}
\lambda_k(v_k, w_k)\to (v,w).
\end{equation*}
Choose $z_k\in \Pi(y_k; \mathcal{F}(x_k))$. By Proposition \ref{l1}, $(v_k, w_k)\in \Hat N((x_k, z_k);\gph \mathcal{F})$ and $\|w_k\|\leq 1$. The cone property of the Fr\'echet normal cone implies
\begin{equation*}
(\lambda_kv_k, \lambda_kw_k)\in \Hat N((x_k, z_k);\gph \mathcal{F}).
\end{equation*}
Since $D(x_k, y_k)=\|y_k-z_k\|\to D(\ox,\oy)=0$, one has that $z_k\to \oy$ as $k\to \infty$. Taking into account the fact that $\lambda_kv_k\to v$ and $\lambda_kw_k\to 0=w$, we obtain $v\in D^*\mathcal{F}(\ox,\oy)(0)$. The inclusion $\subseteq$ has been proved.

Under the Lipschitz continuity of $D$, fix $v\in D^*\mathcal{F}(\ox,\oy)(0)$. Then there exist sequences $(x_k, y_k)\xrightarrow{\gph \mathcal{F}}  (\ox,\oy)$, $(v_k, w_k)\to (v,0)$ with
\begin{equation*}
(v_k, w_k)\in \Hat N((x_k, y_k); \gph \mathcal{F}).
\end{equation*}
By Proposition \ref{l3},  $(v_k, w_k)\in \lambda_k\Hat\partial D(x_k, y_k)$ for $\lambda_k:=\|w_k\|\dn 0$. Thus, $(v, 0)\in \partial^\infty D(\ox,\oy)$. Therefore, the inclusion holds as equality. $\h$ \vspace*{0.05in}

In the proof of (i)$\Rightarrow$ (ii) below, we provide an alternative simplified proof of the Mordukhovich criteria for the Aubin property of set-valued mapping. Similar to the proof given in \cite[Theorem 9.40]{rw}, Lipschitz continuity of the function $D$ from (\ref{scalar}) is employed. However, our proof is based solely on simple analysis from the previous propositions.

\begin{Theorem}\label{E} The following are equivalent:\\[1ex]
{\rm (i)} For any $s\in \Bbb N$ and for any lower semicontinuous function $\psi: \R^s\to (-\infty,\infty]$, one has $\partial^\infty\psi(\ox)=\{0\}$, where $\ox\in \mbox{\rm dom }\psi$, if and only if $\psi$ is Lipschitz continuous around $\ox$.\\
{\rm (ii)} For any $m,n\in \N$, and for any closed-graph set-valued mapping $\mathcal{F}: \R^m\tto\R^n$, one has $D^*\mathcal{F}(\ox,\oy)(0)=\{0\}$, where $(\ox,\oy)\in \gph \mathcal{F}$, if and only if $\mathcal{F}$ is Lipschitz-like around $(\ox,\oy)$.
\end{Theorem}
{\bf Proof.} Let us first assume that (i) is satisfied. Fix $m,n\in \N$ and a closed-graph set-valued mapping $\mathcal{F}: \R^m\tto\R^n$ with $D^*\mathcal{F}(\ox,\oy)(0)=\{0\}$, where $(\ox,\oy)\in \gph \mathcal{F}$. We will show that $\mathcal{F}$ is Lipschitz-like around $(\ox,\oy)$. By Proposition \ref{cg}, the function $D$ defined by $\mathcal{F}$ in (\ref{scalar}) is lower semicontinuous. Using the upper estimate (\ref{upper estimate}), one sees that $\partial^\infty D(\ox,\oy)=\{0\}$ (it always contains $0$). Thus, by (i), $D$ is Lipschitz continuous around $(\ox,\oy)$, which implies that $\mathcal{F}$ is Lipschitz-like around $(\ox,\oy)$ by Theorem \ref{scarl}. Conversely, suppose that $\mathcal{F}$ is Lipschitz-like around $(\ox,\oy)$. Again, by Theorem \ref{scarl}, the function $D$ is Lipschitz continuous around $(\ox,\oy)$. Thus, (\ref{upper estimate}) holds as equality. Moreover, by the converse of (i), $\partial^\infty D(\ox,\oy)=\{0\}$. Therefore, $D^*\mathcal{F}(\ox,\oy)(0)=\{0\}$. The statement (ii) has been proved.

Conversely, suppose that (ii) is satisfied. Fix a lower semicontinuous function $\psi: \R^s\to (-\infty,\infty]$ and define $\mathcal{F}(x)=[\psi(x), \infty)$. By the definition, for $\oy:=\psi(\ox)$, where $\ox\in\mbox{\rm dom }\psi$,
\begin{equation*}\label{function and set}
\partial^\infty\psi(\ox)=D^*\mathcal{F}(\ox, \oy)(0).
\end{equation*}
It is an easy exercise to show that $\psi$ is Lipschitz-like around $\ox$ if and only if $\mathcal{F}$ is Lipschitz-like around $(\ox,\oy)$. By (ii), $\psi$ is Lipschitz continuous around $\ox$ if and only if $D^*\mathcal{F}(\ox, \oy)(0)=\{0\}$, or equivalently $\partial^\infty\psi(\ox)=\{0\}$. $\h$ \vspace*{0.05in}

Let us close the section with simple specification for set-valued mappings with convex graphs.

\begin{Proposition}\label{ls1} Let $\mathcal{F}: \R^m\tto \R^n$ be a convex set-valued mapping. Then the function $D$ defined by {\rm (\ref{scalar})} is convex and
\begin{equation*}
D^*\mathcal{F}(\ox,\oy)(0)=N(\ox, \mbox{\rm dom }\mathcal{F}).
\end{equation*}
\end{Proposition}
\noindent {\bf Proof.} It is not hard to show that $D$ is a convex function. Fix any $x^*\in D^*\mathcal{F}(\ox,\oy)(0)$. Then $(x^*,0)\in N((\ox,\oy);\gph \mathcal{F})$. Thus,
\begin{equation*}
\la x^*, x-\ox\ra +\la 0, y-\oy\ra\leq 0 \mbox{ for all }(x,y)\in \gph \mathcal{F}.
\end{equation*}
For any $x\in \mbox{\rm dom }\mathcal{F}$. Choose $y\in \mathcal{F}(x)$. Then
\begin{equation*}
\la x^*, x-\ox\ra =\la x^*, x-\ox\ra +\la 0, y-\oy\ra\leq 0.
\end{equation*}
Thus, $x^*\in N(\ox;\mbox{\rm dom }\mathcal{F})$. The proof of the converse is also straightforward. $\h$

\begin{Theorem}\label{Lip1} Let $\mathcal{F}: \R^m\tto \R^n$ be a convex set-valued mapping with closed graph. Then the following are equivalent:\\[1ex]
{\rm (i)} $\mathcal{F}$ has the Aubin property around $(\ox,\oy)\in \gph \mathcal{F}$.\\
{\rm (ii)} The function $D$ given by {\rm (\ref{scalar})} is Lipschitz continuous around $(\ox,\oy)$.\\
{\rm (iii)} $\ox\in \mbox{\rm int }(\mbox{\rm dom }\mathcal{F})$.\\
{\rm (iv)} $D^*\mathcal{F}(\ox,\oy)(0)=\{0\}.$
\end{Theorem}
\noindent {\bf Proof.} The equivalence of (i) and (ii) has been stated in Theorem \ref{scarl}. Let us prove the equivalence of (ii) and (iii). We have that
\begin{align*}
\mbox{\rm dom }D&=\{(x,y)\in\R^m\times \R^n\; |\; D(x,y)<\infty\}\\
&=\{(x,y)\in \R^m\times\R^n\; |\; \mathcal{F}(x)\neq \emptyset\}\\
&=\mbox{\rm dom }\mathcal{F}\times \R^n.
\end{align*}
It follows that
\begin{equation*}
\partial^\infty D(\ox,\oy)=N((\ox,\oy); \mbox{\rm dom }D)=N(\ox; \mbox{\rm dom }\mathcal{F})\times \{0\}.
\end{equation*}
Thus, (ii) is equivalent to the fact that $N(\ox; \mbox{\rm dom }\mathcal{F})\times \{0\}=\{(0,0)\}$ or $\ox\in \mbox{\rm int }(\mbox{\rm dom }\mathcal{F})$. The equivalence of (iii) and (iv) follows from Proposition \ref{ls1}. $\h$

\begin{Corollary} Let $f: \R^n\to (-\infty, \infty]$ be a lower semicontinuous convex function and let $\ox\in\mbox{\rm dom }f$. The following are equivalent:\\[1ex]
{\rm (i)} $f$ is Lipschitz continuous around $\ox$.\\
{\rm (ii)} $\ox\in \mbox{\rm int }(\mbox{\rm dom }f)$.\\
{\rm (iii)} $\partial^\infty f(\ox)=\{0\}.$
\end{Corollary}
{\bf Proof. }The results follows from Theorem \ref{Lip1} using $\mathcal{F}(x)=[f(x), \infty)$. Then $\gph \mathcal{F}=\epi f$, $\mbox{\rm dom }\mathcal{F}=\mbox{\rm dom }f$, and $D^*\mathcal{F}(\ox, \oy)(0)=\partial^\infty f(\ox)$, where $\oy:=f(\ox)$. $\h$

\section{Lipschitz Continuity of Minimal Time Functions}

In this section, we are going to provide some examples showing that subdifferential and coderivative criteria presented in the previous section are effective tools for recognizing Lipschitz properties.

Given a nonempty closed bounded convex set $F$, the \emph{Minkowski function} associated with $F$ is given by
\begin{equation*}
\rho_F(x):=\inf\{t\geq 0\; |\; x\in tF\}.
\end{equation*}
It is an easy exercise to show that $\rho_F$ is positively homogeneous and subadditive, and hence convex. Moreover, $\rho_F(0)=0$. Since we do not require that $0\in\mbox{\rm int }F$, it is clear $\rho_F$ is an extended real-valued function. It is an easy exercise to show that $\rho_F$ is lower semicontinuous.

Define
$$F^+:=\{v\in \R^n\; |\; \la v, x\ra \geq 0\mbox{ for all }x\in F\}\; \mbox{\rm and } F^{-}:=-F^+.$$

\begin{Proposition} Consider the  Minkowski function $\rho_F$ with $\ox\in \mbox{\rm dom }\rho_F$. Then
\begin{equation*}
\partial^\infty\rho_F(\ox)=F^{-}\cap\{\ox\}^+.
\end{equation*}
Consequently, $\rho_F$ is Lipschitz continuous around $\ox$ if and only if $F^{-}\cap\{\ox\}^+=\{0\}$. Moreover, the following are equivalent:\\[1ex]
{\rm (i)} $0\in \mbox{\rm int }F$.\\
{\rm (ii)} $F^{-}=\{0\}$.\\
{\rm (iii)} $\rho_F$ is Lipschitz continuous around $0$.\\
{\rm (iv)} $\rho_F$ is globally Lipschitz.
\end{Proposition}
{\bf Proof.} Let us compute $\partial^\infty\rho_F(\ox)$. Since $\rho_F$ is convex and lower semicontinuous, by \cite[Proposition 8.12]{rw}, one has
\begin{equation*}
\partial^\infty\rho_F(\ox)=N(\ox; \mbox{\rm dom }\rho_F),
\end{equation*}
where $\mbox{\rm dom }\rho_F=\mbox{\rm cone }F=\cup_{t\geq 0}tF.$ By the definition, $v\in N(\ox; \mbox{\rm dom }\rho_F)$ if and only if
\begin{equation}\label{cone condition}
\la v, y-\ox\ra\leq 0\;\mbox{\rm for all }y\in \mbox{\rm dom }\rho_F.
\end{equation}
Thus, for any $x\in F$ and $t\geq 0$, one has
\begin{equation*}
\la v, tx-\ox\ra\leq 0.
\end{equation*}
This implies $v\in F^{-}$ by letting $t\to \infty$, and $v\in \{\ox\}^{+}$ by letting $t=0$.

Take any $v\in F^{-}\cap \{\ox\}^+$. It is obvious that $\la v, y\ra \leq 0\leq \la v, \ox\ra$ for any $y\in \mbox{\rm dom }\rho_F$. Thus, (\ref{cone condition}) is satisfied, and hence $v\in N(\ox; \mbox{\rm dom }\rho_F)=\partial^\infty\rho_F(\ox)$. It follows that
\begin{equation*}
\partial^\infty\rho_F(\ox)=F^{-}\cap \{\ox\}^+.
\end{equation*}
Therefore, $\rho_F$ is Lipschitz continuous around $\ox$ if and only if $F^{-}\cap\{\ox\}^+=\{0\}$.

In the case where $\ox=0$, one has $\{\ox\}^+=\R^n$, and hence $\partial^\infty\rho_F(\ox)=F^{-}\cap \{\ox\}^+=F^{-}$. Thus, (ii) and (iii) are equivalent. The fact that (i) and (ii) are equivalent follows from the convex separation theorem, and the fact that (iii) and (iv) are equivalent follows from the fact that $\rho_F$ is positively homogeneous and subadditive with $\rho_F(0)=0$. $\h$ \vspace*{0.05in}

Let us consider another class of functions called the minimal time function that plays an important role in optimization.  Given a \emph{nonempty closed bounded convex set} $F$ (this is our standing assumption on $F$ in this section unless otherwise stated), and given a nonempty closed set $\O$, the minimal time function with dynamic $F$ to $\O$ is given by
\begin{equation}\label{MT}
\mathcal{T}_F(x; \O):=\inf\{t\geq 0\; |\; (x+tF)\cap \O\neq \emptyset\}.
\end{equation}
This function becomes the distance function to $\O$ when $F$ is the closed unit ball of $\R^n$. However, without the assumption that $0\in \mbox{\rm int }F$, this function share less common properties with the distance function. For instance, $\mathcal{T}_F(x; \O)$ is an extended-real valued function and not Lipschitz continuous in general. The readers are referred to \cite{cowo,MN} and the references therein for generalized differentiation properties of the function in infinite dimensions.

We will use a generalized differentiation approach based on singular subgradients to study the Lipschitz continuity of this class of functions. This study continues our recent developments in \cite{NZ}, where a special case of function (\ref{MT}) with $F$ being a singleton has been addressed. Given a set-valued mapping $\mathcal{G}: \R^m\tto \R^n$, let us consider a more general function given by
\begin{equation}\label{MT1}
\mathcal{T}^F_{\mathcal{G}}(x;y):=\inf\{t\geq 0\; |\; (y+tF)\cap \mathcal{G}(x)\neq\emptyset\}.
\end{equation}
For any $(\ox, \oy)\in \mbox{\rm dom }\mathcal{T}_{\mathcal{G}}^F$, define
\begin{equation*}
\Pi^F_{\mathcal{G}}(\ox, \oy)=(\oy+\bar tF)\cap \mathcal{G}(\ox),
\end{equation*}
where $\bar t:=\mathcal{T}_{\mathcal{G}}^F(\ox,\oy)$.

The following proposition can be easily proved following \cite{MN}.
\begin{Proposition} Let $\mathcal{G}: \R^m\tto \R^n$ be a set-valued mapping with closed graph. Then\\[1ex]
{\rm (i)} $\mathcal{T}^F_{\mathcal{G}}(\ox;\oy)=0$ iff $\oy\in \mathcal{G}(\ox)$.\\
{\rm (ii)} $\mathcal{T}^F_{\mathcal{G}}$ is lower semicontinuous.\\
{\rm (iii)} For any $(\ox, \oy)\in \mbox{\rm dom }\mathcal{T}_{\mathcal{G}}^F$, the projection set $\Pi^F_{\mathcal{G}}(\ox, \oy)$ is nonempty.\\
{\rm (iv)} $\mathcal{T}^F_{\mathcal{G}}$ is convex if $\mathcal{G}$ has a convex graph.
\end{Proposition}

\begin{Proposition}\label{MT2} Let $\mathcal{G}: \R^m\tto \R^n$ be a set-valued mapping with closed graph. Consider the minimal time function {\rm (\ref{MT1})}. Suppose that $\gph\mathcal{G}$ is closed and $\oy\notin\mathcal{G}(\ox)$. For any $\oz\in \Pi^F_{\mathcal{G}}(\ox, \oy)$,
\begin{equation*}
\Hat\partial \mathcal{T}_{\mathcal{G}}^F(\ox,\oy)\subseteq \{(u,v)\in \R^m\times\R^n\; |\; (u,v)\in \Hat N((\ox, \oz); \gph \mathcal{G}), v\in S^*\},
\end{equation*}
where $S^*:=\{v\in \R^n\; |\; \sigma_F(-v):=\sup_{x\in F}\la -v,x\ra= 1\}$.
\end{Proposition}
{\bf Proof. }Fix any $(u, v)\in \Hat\partial \mathcal{T}_{\mathcal{G}}^F(\ox,\oy)$. For any $\ve>0$, there exists $\delta>0$ such that
\begin{equation}\label{FS}
\la u, x-\ox\ra+\la v, y-\oy\ra \leq \mathcal{T}_{\mathcal{G}}^F(x,y)-\mathcal{T}_{\mathcal{G}}^F(\ox,\oy)+\ve(\|x-\ox\|+\|y-\oy\|)
\end{equation}
whenever $\|x-\ox\|<\delta$ and $\|y-\oy\|<\delta$. Fix $\oz\in \Pi^F_{\mathcal{G}}(\ox, \oy)$ and $(x, z)\in \gph \mathcal{G}$ with $\|x-\ox\|<\delta$ and $\|z-\oz\|<\delta$. Then $\|(z-\oz+\oy)-\oy\|<\delta$. Thus,
\begin{equation*}
\la u, x-\ox\ra +\la v, z-\oz\ra\leq \mathcal{T}_{\mathcal{G}}^F(x, z-\oz+\oy)-\mathcal{T}_{\mathcal{G}}^F(\ox,\oy)+\ve(\|x-\ox\|+\|z-\oz\|).
\end{equation*}
Since $\oz\in \oy+\bar tF$, where $\bar t:=\mathcal{T}_{\mathcal{G}}^F(\ox,\oy)$, one has $z\in (z-\oz+\oy+\bar tF)\cap \mathcal{G}(x)$. Thus, $\mathcal{T}_{\mathcal{G}}^F(x, z-\oz+\oy)\leq \bar t$. It follows that
\begin{equation*}
\la u, x-\ox\ra +\la v, z-\oz\ra\leq \ve(\|x-\ox\|+\|z-\oz\|).
\end{equation*}
Therefore, $(u, v)\in \Hat N((\ox, \oz); \gph \mathcal{G})$.

Fix $x=\ox$ in (\ref{FS}) and let $\O:=\mathcal{G}(\ox)$. One has that $v\in \Hat\partial \mathcal{T}_\O^F(\ox)$. Then $\sigma_F(-v)= 1$ by \cite{MN}. $\h$

\begin{Theorem}\label{LC} Let $\mathcal{G}: \R^m\tto \R^n$ be a set-valued mapping with closed graph. Suppose that $(\ox, \oy)\in \gph \mathcal{G}$. Assume that the following implication holds:
\begin{equation}\label{im}
u\in D^*\mathcal{G}(\ox, \oy)(v) \; \mbox{\rm and }v\in F^{-}\Rightarrow [u=v=0].
\end{equation}
 Then the minimal time function {\rm (\ref{MT1})} is Lipschitz continuous around $(\ox,\oy)$. The converse also holds true if we assume additionally that $0\in F$.
\end{Theorem}
{\bf Proof. }Suppose that the implication (\ref{im}) holds. Fix any $(u,v)\in \partial^\infty \mathcal{T}^F_\mathcal{G}(\ox, \oy)$. Then there exist sequences $\lambda_k\dn 0$, $(u_k, v_k)\xrightarrow {\mathcal{T}^F_\mathcal{G}}(\ox, \oy)$, $(u_k, v_k)\in \Hat\partial \mathcal{T}^F_\mathcal{G}(x_k,y_k)$ with $\lambda_k(u_k, v_k)\to (u,v)$. In the case where $(x_k, y_k)\in \gph \mathcal{G}$ for a subsequence (without relabeling), by \cite[Theorem 3.1 (a)]{JH},
\begin{equation*}
(u_k, v_k)\in \Hat N((x_k, y_k); \gph \mathcal{G})\; \mbox{ and }\sigma_F(-v_k)\leq 1.
\end{equation*}
It follows that $\lambda_k(u_k, v_k)\in \Hat N((x_k, y_k); \gph \mathcal{G})$ and $\sigma_F(-\lambda_k v_k)\leq \lambda_k$ since $\sigma_F(\cdot)$ is positively homogeneous. Since the function $\sigma_F(\cdot)$ is lower semicontinuous, one has
\begin{equation*}
(u, v)\in N((\ox, \oy); \gph \mathcal{G})\; \mbox{ and }\sigma_F(-v)\leq 0.
\end{equation*}
It follows that $u\in D^*\mathcal{G}(\ox, \oy)(-v)$ and $-v\in F^{-}$. Thus, $(u, v)=(0,0)$. Therefore, $\partial^\infty\mathcal{T}^F_{\mathcal{G}}(\ox, \oy)=\{(0,0)\}$, and hence $\mathcal{T}^F_{\mathcal{G}}$ is Lipschitz continuous around $(\ox,\oy)$. In the case where $(x_k, y_k)\notin\gph \mathcal{G}$, we use Proposition \ref{MT2} and arrive at the same conclusion.

Let us prove the converse. Suppose that $\mathcal{T}^F_{\mathcal{G}}$ is Lipschitz continuous around $(\ox, \oy)$. Let $(u, v)$ satisfy (\ref{im}). Then
\begin{equation*}
(u, -v)\in N((\ox, \oy); \gph \mathcal{G}),
\end{equation*}
for some $v\in F^{-}$. There exist sequences $(x_k, y_k)\xrightarrow{\gph \mathcal{G}}(\ox, \oy)$, $(u_k, -v_k)\rightarrow (u, -v)$ and $(u_k, -v_k)\in \Hat N((x_k, y_k); \gph \mathcal{G})$. Since $v_k\to v$ and $0\in F$, by extracting a subsequence, we can assume that $0\leq \lambda_k:=\sigma_F(v_k)\to 0$ as $k\to \infty$. Indeed, using the compactness of $F$, for each $k$, choose $f_k\in F$ such that $\lambda_k=\la v_k, f_k\ra$. By extracting a subsequence, we can assume $(f_k)$ converges to some $f\in F$, and hence $\lambda_k=\la v_k, f\ra\to \la v, f\ra$. Since $f\in F^{-}$, $\la v, f\ra \leq 0$. Moreover, $\lambda_k\geq 0$ for every $k$. Thus, $\lambda_k\to 0$. In the case where $\lambda_k>0$ for a subsequence (without relabeling), one has $\sigma_F(\dfrac{v_k}{\lambda_k})=1$ and
\begin{equation}\label{NC}
(\dfrac{u_k}{\lambda_k}, \dfrac{-v_k}{\lambda_k})\in \Hat N((x_k, y_k); \gph \mathcal{G}).
\end{equation}
By \cite[Theorem 3.1 (b)]{JH}, $(\dfrac{u_k}{\lambda_k}, \dfrac{-v_k}{\lambda_k})\in \Hat\partial \mathcal{T}^F_\mathcal{G}(x_k, y_k)$, or $(u_k, -v_k)\in \lambda_k\Hat\partial \mathcal{T}^F_\mathcal{G}(x_k, y_k)$. This implies $(u, -v)\in \partial^\infty T^F_\mathcal{G}(\ox,\oy)=\{(0,0)\}$, and hence $u=v=0$.

In the case where $\sigma_F(v_k)=0$ for $k$ sufficiently large. Choose a sequence of positive numbers $\lambda_k\dn 0$. Then we also have $\sigma_F(\dfrac{v_k}{\lambda_k})=0\leq 1$, and (\ref{NC}) also holds. Using the same argument, we also arrive at $u=v=0$. The proof is now complete. $\h$

\begin{Corollary} Let $\mathcal{G}: \R^m\tto \R^n$ be a set-valued mapping with closed graph. Fix $(\ox, \oy)\in \gph \mathcal{G}$. Suppose that $0\in \mbox{\rm int }F$. Then the minimal time function {\rm(\ref{MT1})} is Lipschitz continuous around $(\ox, \oy)$ if and only if $\mathcal{G}$ is Lipschitz-like around $(\ox,\oy)$.
\end{Corollary}
{\bf Proof.} Under the assumption $0\in \mbox{\rm int }F$, one has $F^{-}=\{0\}$. Thus, condition (\ref{im}) requires $D^*\mathcal{G}(\ox, \oy)(0)=\{0\}$, which is the necessary and sufficient condition for $\mathcal{G}$ to be Lipschitz-like around $(\ox, \oy)$. $\h$ \vspace*{0.05in}

Let us derive a necessary and sufficient condition for Lipschitz continuity of the minimal time function (\ref{MT}).
\begin{Corollary} Let $\O$ be a nonempty closed set and let $\oy\in \O$. If $N(\oy; \O)\cap F^{+}=\{0\}$, then the minimal time function {\rm(\ref{MT})} is Lipschitz continuous around $\oy$. The converse also holds true if $0\in F$.
\end{Corollary}
{\bf Proof. }Consider the set-valued mapping defined by $\mathcal{G}(x)=\O$. Then $\gph \mathcal{G}=\R^m\times \O$. Fix any $\ox\in \R^m$. Then $N((\ox, \oy); \gph \mathcal{G})=\{0\}\times N(\oy; \O)$. Moreover,
\begin{equation*}
\mathcal{T}_\mathcal{G}^F(x,y)=\mathcal{T}^F_\O(y) \; \mbox{\rm for all }(x,y)\in \R^m\times \R^n.
\end{equation*}
Suppose that $N(\oy; \O)\cap F^{+}=\{0\}$. Let us show that condition (\ref{im}) is satisfied. Fix $u\in D^*\mathcal{G}(\ox, \oy)(v) \; \mbox{\rm and }v\in F^{-}$. Then $-v\in F^{+}$ and $(u, -v)\in N((\ox,\oy); \gph \mathcal{G})$. This implies $u=0$ and $-v\in N(\oy; \O)\cap F^{+}$, which implies $u=v=0$. Therefore, $\mathcal{T}^F_\O$ is Lipschitz continuous around $\oy$. The proof of the converse also follows from Theorem \ref{LC} by a similar argument. $\h$ \vspace*{0.05in}

Let us now establish sufficient conditions that guarantee the Lipschitz continuity of the minimal time function (\ref{MT1}) for $(\ox, \oy)\notin\gph \mathcal{\mathcal{G}}$.

\begin{Theorem}\label{LC2} Let $\mathcal{G}: \R^m\tto \R^n$ be a set-valued mapping with closed graph. Let $\oy\notin \mathcal{\mathcal{G}}(\ox)$ with $(\ox,\oy)\in \mbox{\rm dom }\mathcal{T}^F_{\mathcal{G}}$. Suppose that for any $\oz\in \Pi^F_{\mathcal{G}}(\ox, \oy)$, the following implication holds:
\begin{equation}\label{im1}
u\in D^*\mathcal{G}(\ox, \oz)(v) \; \mbox{\rm and }v\in F^{\perp}\Rightarrow [u=v=0],
\end{equation}
where $F^{\perp}:=\{v\in \R^n\; |\; \la v, x\ra=0\;\mbox{\rm for all }x\in F\}$. Then the minimal time function {\rm (\ref{MT1})} is Lipschitz continuous around $(\ox,\oy)$.
\end{Theorem}
{\bf Proof. }Let us show that $\partial^\infty\mathcal{T}^F_{\mathcal{G}}(\ox, \oy)=\{(0,0)\}$ under the assumption made. Fix any $(u, v)\in \partial^\infty\mathcal{T}^F_{\mathcal{G}}(\ox, \oy)$, find $\lambda_k\dn 0$, $(x_k, y_k)\xrightarrow{\mathcal{T}^F_{\mathcal{G}}}(\ox,\oy)$, $(u_k, v_k)\in \Hat \partial \mathcal{T}^F_{\mathcal{G}}(x_k, y_k)$ such that
\begin{equation*}
\lambda_k(u_k, v_k)\rightarrow (u,v).
\end{equation*}
Fix $z_k\in \Pi^F_{\mathcal{G}}(x_k, y_k)$ and let $t_k:=\mathcal{T}^F_{\mathcal{G}}(x_k, y_k)$. Then $t_k\to \bar t:=\mathcal{T}^F_{\mathcal{G}}(\ox,\oy)$. Since $z_k\in (y_k+t_kF)\cap \mathcal{G}(x_k)$, it is clear that $(z_k)$ is bounded and one can assume, without loss of generality, that $z_k\to \oz$ as $k\to \infty$. Thus, $\oz\in (\oy+\bar tF)\cap \mathcal{G}(\ox)=\Pi^F_{\mathcal{G}}(\ox, \oy)$ since $\gph \mathcal{G}$ is closed. By Proposition \ref{MT2} and \cite[Theorem 3.2]{JH},
\begin{equation*}
(u_k, v_k)\in \Hat N((x_k, z_k); \gph \mathcal{G})\; \mbox{\rm and }\sigma_F(-v_k)= 1.
\end{equation*}
This implies
\begin{equation*}
(\lambda_k u_k, \lambda_kv_k)\in \Hat N((x_k, z_k); \gph \mathcal{G})\; \mbox{\rm and }v\in F^{\perp}.
\end{equation*}
Thus,
\begin{equation*}
(u, v)\in  N((\ox, \oz); \gph \mathcal{G})\; \mbox{\rm and } v\in F^{\perp}.
\end{equation*}
Then $u\in D^*\mathcal{G}(\ox, \oz)(-v)$ and $-v\in F^{\perp}$, and hence $(u,v)=\{(0,0)\}$ by (\ref{im1}). Therefore, $\partial^\infty\mathcal{T}^F_{\mathcal{G}}(\ox, \oy)=\{(0,0)\}$, which implies the Lipschitz continuity of $\mathcal{T}^F_{\mathcal{G}}$ at $(\ox, \oy)$. $\h$ \vspace*{0.05in}

The following corollaries follow directly from Theorem \ref{LC2}.
\begin{Corollary} Let $\mathcal{G}: \R^m\tto \R^n$ be a set-valued mapping with closed graph. Let $\oy\notin \mathcal{G}(\ox)$ with $(\ox,\oy)\in \mbox{\rm dom }\mathcal{T}^F_{\mathcal{G}}$. Suppose that $0\in \mbox{\rm int }F$ and for any $\oz\in \Pi^F_{\mathcal{G}}(\ox, \oy)$, the set-valued mapping $\mathcal{G}$ is Lipschitz-like around $(\ox, \oz)$. Then the minimal time function {\rm (\ref{MT1})} is Lipschitz continuous at $(\ox,\oy)$.
\end{Corollary}

\begin{Corollary}\label{lips4} Let $\O$ be a nonempty closed subset of $\R^n$. Consider the minimal time function {\rm (\ref{MT})}. Let $\oy\notin \O$ with $\bar t:=\mathcal{T}_\O^F(\oy)<\infty$. If $N(\oz; \O)\cap F^{\perp}=\{0\}$ for all $\oz\in \Pi^F_\O(\oy; \O):=(\oy+\bar tF)\cap \O$, then the minimal time function {\rm (\ref{MT})} is Lipschitz continuous around $\oy$.
\end{Corollary}

Observe that the opposite implication in Theorem \ref{LC2} does not hold true in general as shown in the example below.
\begin{Example}{\rm In $\R^2$, let $\O=\{(x, 0)\; |\; x\in \R\}$, $F=\{0\}\times [0, 1]$, and $\oy=(0, -2)$. Then $\Pi^F_\O(\oy; \O)=\{(0,0)\}$. It is not hard to see that $\mathcal{T}^F_\O$ is Lipschitz continuous at $\oy$, but for $\oz:=(0,0)$, one has
\begin{equation*}
N(\oz; \O)\cap F^{\perp}=\R\times\{0\}.
\end{equation*}}
\end{Example}

For a given vector $d\in \R^n$, $d\neq 0$, and a nonempty closed set $\O\subseteq \R^n$, the \emph{scalarization function} defined by the direction $d$ and target set $\O$ is defined by
\begin{equation}\label{s}
\ph_d(x;\O):=\inf\{t\in\Bbb R\; |\; x+td\in \O\}.
\end{equation}
The only difference in this definition compared with the corresponding minimal time function (\ref{MT}) in which $F=\{d\}$ is that $t$ can take negative values. The scalarization function (\ref{s}) has been extensively used in vector optimization; see, e.g., \cite{TZ} and the references therein.

Following \cite{TZ}, we say that $\O$ satisfies the free-disposal condition in the direction $d$ if $$\O+\R^{+}d=\O.$$
We also say that $\O$ is \emph{normal regular} at $\ox\in \O$ if $\Hat N(\ox; \O)=N(\ox; \O)$. This assumption is automatically satisfied when $\O$ is convex.

We will show that the  property in the proposition below holds true without assuming the free-disposal condition. This improves the corresponding result from \cite{TZ}.
\begin{Proposition}\label{pr} Consider the scalarization function {\rm (\ref{s})}. For any $x\in \R^n$, one has
\begin{equation*}
\ph_d(x+\alpha d; \O)=\ph_d(x; \O)-\alpha.
\end{equation*}
\end{Proposition}
{\bf Proof. }Suppose $|\ph_d(x;\O)|<\infty$. Then $x+\ph_d(x; \O)d\in \O$. Thus $x+\alpha d+(\ph_d(x; \O)-\alpha)d\in \O$, so $\ph_d(x+\alpha d; \O)\leq \ph_d(x; \O)-\alpha<\infty$.
Let us show that $\ph_d(x+\alpha d; \O)>-\infty$. If so, let $t_k\to -\infty$ and $x+\alpha d +t_kd\in \O$. Then $\ph_d(x; \O)\leq \alpha +t_k$, which is a contradiction. Then
$$\ph_d(x; \O)=\ph_d(x+\alpha d-\alpha d; \O)\leq \ph_d(x+\alpha d; \O)+\alpha.$$
The proof is complete for this case. Suppose $\ph_d(x; \O)=\infty$. Then $x\notin \O-\Bbb R d$. Thus $x+\alpha d\notin \O-\Bbb R d$. So $\ph_d(x+\alpha d; \O)=\infty$. Now suppose $\ph_d(x; \O)=-\infty$. Then let $t_k\to -\infty$ and $x+t_kd\in \O$. Then $x+\alpha d+(t_k-\alpha)d\in \O$, so $\ph_d(x+\alpha d; \O)\leq t_k-\alpha$. This implies $\ph_d(x+\alpha d; \O)=-\infty$. $\h$ \vspace*{0.05in}

Let us know obtain representations of the Fr\'echet and limiting subdifferential of the scalarization function (\ref{s}) without assuming the free-disposal condition.

\begin{Proposition}\label{sub calculation} Let $\ox\in \mbox{\rm dom }\ph_d(\cdot; \O)$ and $\tilde{x}:=\ox+td$, where $t:=\ph_d(\ox; \O)$. Then
\begin{equation}\label{fr}
\Hat\partial \ph_d(\ox; \O)\subseteq\{w\in \R^n\; |\; \la w, -d\ra =1\}\cap \Hat N(\tilde{x};\O),
\end{equation}
and
\begin{equation}\label{lm}
\partial\ph_d(\ox; \O)\subseteq\{w\in \R^n\; |\; \la w, -d\ra =1\}\cap N(\tilde{x};\O),
\end{equation}
The inclusion {\rm (\ref{fr})} holds as equality if we assume additionally that $\ph_d(\cdot; \O)$ is calm at $\ox$, and the inclusion {\rm (\ref{lm})} holds as equality if $\ph_d(\cdot; \O)$ is calm at $\ox$ and $\O$ is normal regular at $\tilde{x}$.
\end{Proposition}
{\bf Proof.} Fix any $w\in\Hat\partial \ph_d(\ox; \O)$. For any $\ve>0$, there exists $\delta>0$ such that
\begin{equation*}
\la w, x-\ox\ra \leq \ph_d(x; \O)-\ph_d(\ox; \O) +\ve \|x-\ox\|\mbox{ whenever }x\in \B(\ox; \delta).
\end{equation*}
For $t>0$ sufficiently small, one has $x:=\ox+td\in \B(\ox; \delta)$, and hence
\begin{align*}
t\la w, d\ra=\la w, \ox+td-\ox\ra &\leq \ph_d(\ox+td; \O)-\ph_d(\ox; \O)+\ve t \|d\|\\
&=\ph_d(\ox; \O) -t-\ph_d(\ox; \O)+\ve t\|d\|=-t+\ve t\|d\|.
\end{align*}
This implies $\la w, -d\ra\geq  1$. Using $x-tv$ instead of $x+tv$ in a similar one, one obtain $\la w, -d\ra\leq 1$, and hence $\la w, -d\ra=1$.

For any $x\in \O$, one has $\ph_d(x; \O)\leq 0$. Fix any  $x\in \B(\tilde{x}; \delta)\cap \O$. Then $\|(x-\tilde{x}+\ox)-\ox\|<\delta$. So
\begin{align*}
\la w, x-\tilde{x}\ra&=\ph_d(x-\tilde{x}+\ox; \O)-\ph_d(\ox; \O)+\ve \|x-\tilde{x}\|\\
&=\ph_d(x-td; \O)-\ph_d(\ox; \O)+\ve \|x-\tilde{x}\|\\
&=\ph_d(x; \O)-t+\ph_d(\ox; \O)+\ve \|x-\tilde{x}\|\\
&\leq \ve \|x-\tilde{x}\|.
\end{align*}
This implies $w\in \Hat N(\tilde{x}; \O)$. The inclusion $\subseteq$ has been proved.

The first inclusion has been proved. Let us prove the opposite inclusion. Fix any $w\in \R^n$ with $\la w, -d\ra=1$ and $w\in N(\tilde{x}; \O)$. For any $\ve>0$, there exists $\delta>0$ and a constant $\ell>0$ such that
\begin{equation*}
\la w, x-\tilde{x}\ra \leq \ve \|x-\tilde{x}\|\; \mbox{\rm whenever }x\in \B(\tilde{x}; \delta)\cap \O,
\end{equation*}
and  $\|\ph_d(x; \O)-\ph_d(\ox; \O)\|\leq \ell\|x-\ox\|$ whenever $\|x-\ox\|<\delta$. Let $\gamma:=\dfrac{\delta}{1+\ell\|d\|}$. If $\|x-\ox\|<\gamma$, then
$$\|x+\ph_d(x; \O)d-\tilde{x}\|=\|x+\ph_d(x; \O)d-(\ox+\ph_d(\ox; \O)d)\|\leq \|x-\ox\|+\ell\|x-\ox\|\|d\|< \delta.$$
Since $x+\ph_d(x; \O)d\in \O$ as always, it follows that
\begin{equation*}
\la w,  x+\ph_d(x; \O)d-\tilde{x}\ra \leq \ve \|x+\ph_d(x; \O)d-\tilde{x}\|\leq \ve (1+\ell\|d\|)\|x-\ox\|.
\end{equation*}
Thus,
\begin{align*}
\la w,  x-\ox\ra &\leq (\ph_d(x; \O)-\ph_d(\ox; \O))\la w, -d\ra+ \ve (1+\ell\|d\|)\|x-\ox\|\\
&=\ph_d(x; \O)-\ph_d(\ox; \O)+ \ve (1+\ell\|d\|)\|x-\ox\|.
\end{align*}
Since $\ve$ is arbitrary, we have that $w\in \Hat\partial \ph_d(\ox; \O)$.

The proof for the inclusion $\subseteq$ in the limiting subdifferential representation (\ref{lm}) follows from a simple limiting procedure. Since $\O$ is regular at $\tilde{x}$, applying the definition and  (\ref{fr}), one has
\begin{align*}
\{w\in \R^n\; |\; \la w, -d\ra =1\}\cap N(\tilde{x};\O)&=\{w\in \R^n\; |\; \la w, -d\ra =1\}\cap N(\tilde{x};\O)=\Hat\partial\ph_d(\ox; \O).
\end{align*}
Because $\Hat\partial\ph_d(\ox; \O)\subseteq \partial \ph_d(\ox; \O)$, the opposite inclusion of (\ref{lm}) follows. The proof is now complete. $\h$

\begin{Proposition} Consider the minimal time function {\rm (\ref{s})}. Fix $\ox\in \mbox{\rm dom }\ph_d(\cdot; \O)$ and denote $\tilde{x}:=\ox+td$, where $t:=\ph_d(\ox; \O)$. Suppose that $N(\tilde{x}; \O)\cap \{d\}^\perp =\{0\}$. Then the minimal time function {\rm (\ref{s})} is Lipschitz continuous around $\ox$. The converse also holds true if we assume additionally that $\O$ is normal regular at $\tilde{x}$.
\end{Proposition}
{\bf Proof. }Suppose that $N(\tilde{x}; \O)\cap \{d\}^\perp =\{0\}$. Let us show that the minimal time function {\rm (\ref{s})} is Lipschitz continuous around $\ox$.  Fix any $v\in \partial^\infty\ph_d(\ox; \O)$. Then there exist sequences $\lambda_k\dn 0$, $x_k\xrightarrow{\ph_d(\cdot; \O)}\ox$, $v_k\in \Hat\partial \ph_d(x_k; \O)$ with $\lambda_kv_k\to v$. Let $\tilde{x}_k:=x_k+t_kd$, where $t_k:=\ph_d(x_k; \O)$. Then $t_k\to t$ and $\tilde{x}_k\xrightarrow{\O}\ox$. By (\ref{fr}), one has
\begin{equation*}
\la v_k, d\ra =-1\; \mbox{\rm and } v_k\in \Hat N(\tilde{x}_k; \O).
\end{equation*}
This implies
\begin{equation*}
\la \lambda_k v_k, d\ra =-\lambda_k \;\mbox{\rm and } \lambda_k v_k\in \Hat N(\tilde{x}_k; \O).
\end{equation*}
It follows by letting $k\to\infty$ that $\la v, d\ra =0$ and $v\in N(\ox; \O)$, so $v=0$. Therefore, $\partial^\infty\ph_d(\ox; \O)=\{0\}$, and hence $\ph_d(\cdot; \O)$ is Lipschitz continuous around $\ox$ by Theorem \ref{mcriterionf}.

For the converse, we assume that $\ph_d(\cdot; \O)$ is Lipschitz continuous around $\ox$ and $\O$ is regular at $\ox$. Then $\partial \ph_d(\ox; \O)=\Hat\partial \ph_d(\ox; \O)$ under the assumptions made, which guarantee the equality in (\ref{fr}) and (\ref{lm}). It is well-known and not hard to see that $\partial \ph_d(\ox; \O)$ is nonempty under the Lipschitz continuity of $\ph_d(\cdot;\O)$ around $\ox$. Choose $w\in \partial \ph_d(\ox; \O)=\Hat\partial \ph_d(\ox; \O)$. Then $\la w, -d\ra =1$ and $w\in \Hat N(\tilde{x}; \O)$ by Proposition \ref{sub calculation}. Fix any $v\in N(\tilde{x}; \O)\cap \{d\}^\perp=\Hat N(\tilde{x}; \O)\cap \{d\}^\perp$. Since $\Hat N(\tilde{x}; \O)$ is a convex cone,  $w+kv\in \Hat N(\tilde{x}; \O)$, and, moreover, $\la w+kv, -d\ra =\la w, -d\ra +k\la w, -d\ra =1$. It follows that $w+kv\in \Hat\partial \ph_d(\ox)$, so
$$\dfrac{1}{k}w+v\in \dfrac{1}{k}\Hat\partial \ph_d(\ox; \O).$$
Thus, $v\in \partial^\infty \ph_d(\ox; \O)=\{0\}$ because $\ph_d(\cdot; \O)$ is Lipschitz continuous around $\ox$. We have shown that $N(\tilde{x}; \O)\cap \{d\}^\perp =\{0\}$. $\h$


\begin{thebibliography}{99}

\bibitem{af} Aubin, J.-P: Lipschitz behavior of solutions to convex minimization problem, Math. Oper. Res. \textbf{9}, 87--111 (1984).

\bibitem{c-1983} Clarke, F.H.: Optimization and Nonsmooth
Analysis. Wiley, New York (1983).

\bibitem{cowo} Colombo, G., Wolenski, P.R.: The subgradient formula for the minimal time
function in the case of constant dynamics in Hilbert space. J.
Global Optim. \textbf{28}, 269–-282 (2004).

\bibitem{JH} Jiang, Y., He., Y.:  Subdifferential properties for a class of minimal time functions with moving target sets in normed spaces, Appl. Anal. {\bf 3}, 491--502 (2012).

\bibitem{m76} Mordukhovich, B. S.: Maximum principle in problems of
time optimal control with nonsmooth constraints. J. Appl. Math.
Mech. {\bf 40}, 960--969 (1976).

\bibitem{m88} Mordukhovich, B. S.: Approximation Methods in Problems of Optimization and Control. Nauka, Main Physical and Mathematical Editions, Moscow (1988) (Russian).


\bibitem{m92}  Mordukhovich, B. S.: Sensitivity analysis in nonsmooth optimization, in Theoretical Aspects of Industrial Design (D.A. Field and V. Komkov, eds.). SIAM Proc. in Applied Mathematics \textbf{58}, 32--42, SIAM Publications, Philadelphia, PA (1992).


\bibitem{m93} Mordukhovich, B. S.: Complete characterizations of
openness, metric regularity, and Lipschitzian properties of
multifunctions. Trans. Amer. Math. Soc. {\bf 340}, 1--36 (1993).

\bibitem{mor} Mordukhovich, B.S.: Variational Analysis and
Generalized Differentiation, I: Basic Theory, II: Applications,
Grundlehren Series (Fundamental Principles of Mathematical
Sciences), Vols. 330 and 331, Springer, Berlin (2006).
\bibitem{MN} Mordukhovich, B.S., Nam, N.M.: Subgradients of minimal time functions
under minimal assumptions. J. Convex Anal.{\bf 18}, 915-947 (2011).
\bibitem{mn1} Mordukhovich, B. S., Nam, N. M.: Subgradients of
distance functions with applications to Lipschitzian stability.
Math. Program. {\bf 104}, 635--668 (2005).

\bibitem{NZ} Nam, N. M., Zalinescu, C.: Variational analysis of directional minimal time functions and applications to location problems, to appear in Set-Valued and Variational Analysis.


\bibitem{r85} Rockafellar, R. T.: Lipschitzian property of multifunctions. Nonlinear Anal. \textbf{9}, 867--885 (1985).
\bibitem{rw} Rockafellar, R. T., Wets, R. J.-B.: Variational
Analysis. Springer, Berlin (1998).

\bibitem {TZ}Tammer, Chr., Z\u{a}linescu, C.: Lipschitz properties of the
scalarization function and applications, Optimization \textbf{59}, 305--319 (2010).


\end{thebibliography}
\end{document}